\documentclass[11pt]{article}

\usepackage{authblk} 

\usepackage[ngerman, english]{babel} 
\usepackage[utf8]{inputenc} 
\usepackage[T1]{fontenc}
\usepackage{comment}
\usepackage{verbatim} 

\usepackage{lmodern} 

\usepackage[pdfpagelabels, hidelinks]{hyperref} 
\usepackage{enumerate} 

\usepackage{amsmath, amsthm, amssymb,amsfonts} 
\usepackage{dsfont} 
\usepackage{mathtools}

\usepackage[table]{xcolor}
\usepackage{diagbox}

\newcommand\blfootnote[1]{
	\begingroup
	\renewcommand\thefootnote{}\footnote{#1}
	\addtocounter{footnote}{-1}
	\endgroup
}

\usepackage{xcolor} 

\usepackage{orcidlink}

\theoremstyle{plain} 
\newtheorem{theorem}{Theorem}[section]
\newtheorem{lemma}[theorem]{Lemma}
\newtheorem{proposition}[theorem]{Proposition}
\newtheorem{corollary}[theorem]{Corollary}
\newtheorem{fact}[theorem]{Fact}
\newtheorem{observation}[theorem]{Observation}

\newtheorem*{lemma*}{Lemma}
\newtheorem*{theorem*}{Theorem}
\newtheorem*{corollary*}{Corollary}
\newtheorem*{claim*}{Claim}

\theoremstyle{definition}

\newtheorem{remark-notation}[theorem]{Remark and Notation}

\newtheorem{remark}[theorem]{Remark}
\newtheorem{conjecture}[theorem]{Conjecture}
\AtEndEnvironment{rem}{\hfill\ensuremath{\diamond}}
\newtheorem{question}[theorem]{Question}

\newtheorem*{definition*}{Definition}
\newtheorem*{example*}{Example}
\newtheorem*{remark*}{Remark}

\theoremstyle{definition} 
\newtheorem{open-question}{Question}[section]

\numberwithin{equation}{section}

\title{\Large\textbf{{On Tameness, Measurability and the Independence Property}}}

\author[i,ii] {Lothar Sebastian Krapp \orcidlink{0000-0003-3102-1923}\,}
\author[iii]{Matthieu Vermeil \orcidlink{0009-0009-3280-5907}\,}
\author[ii]{Laura Wirth \orcidlink{0000-0003-2871-5676}\,}
\affil[i]{\,Institut für Interdisziplinäre Sprachevolutionswissenschaft, Universität Zürich, Switzerland}
\affil[ii]{\,Fachbereich Mathematik und Statistik, Universität Konstanz, Germany}
\affil[iii]{\,UFR de Mathématiques, Université Paris Cité, France}

\date{}

\usepackage[skip=5pt+2pt-2pt]{parskip}

\begin{document}
	
	\pagenumbering{arabic}
	
	
	\maketitle
	\blfootnote{Math Subject Classification (2020): Primary 03C64, 28A05; Secondary 12L12, 12J15, 03C40, 03C45, 54H05. Keywords: Borel, first-order definable, NIP, non-measurable, tame geometry. French Title: Modération, mesurabilité et propriété d’indépendance.
	}
	\blfootnote{Corresponding Author: Laura Wirth, laura.wirth@uni-konstanz.de.
	}
	\vspace{-1.3cm}
	
	\renewcommand\abstractname{Abstract}
	\begin{abstract}
		\noindent In the area of Tame Geometry, different model-theoretic tameness conditions are established and their relationships are analyzed. We construct a subfield $K$ of the real numbers that lacks several of such tameness properties. As our main result, we present a first-order formula in the language of rings that defines a non-Borel set in $K$. Moreover, $K$ has the independence property and admits both archimedean and non-archimedean orderings.
	\end{abstract}
	
	\tableofcontents

	
	\section{Introduction}
	
	Several notions of tameness within Model Theory originate from the properties of (first-order) definable sets and functions in o-minimal expansions of the real field or in more general o-minimal structures.\footnote{See van~den~Dries~\cite{Van-den-Dries} for a comprehensive treatment of tameness in o-minimal structures. See Hieronymi~\cite{Hieronymi} for a survey on the origins and developments of Tame Geometry for expansions of the reals.
	}
	This research note is motivated by the first and third author's work \cite{Krapp-Wirth}, in which model-theoretic tameness properties are applied in order to re-examine and strengthen the Fundamental Theorem of Statistical Learning in the context of o-minimality (see the main result \cite[Theorem~4.7]{Krapp-Wirth}). This line of research thus aims to deepen the connection between Tame Geometry and the theoretical underpinnings of Machine Learning.
	In this context, the following two observations are crucial:
	\begin{enumerate}
		\item Due to the Cell Decomposition Theorem, any definable set in an o-minimal ordered field is Borel with respect to the order topology (cf.\ Kaiser~\cite[Proposition~1.1]{Kaiser}, Karpinski and Macintyre~\cite[Lemma~6]{Karp-Mac}). This ensures that several essential measurability conditions, as introduced in \cite[Section~3]{Krapp-Wirth}, are readily satisfied by definable sets in o-minimal expansions of the reals. Measurability, in turn, allows certain probabilities to be evaluated, which is key in the concept of \emph{probably approximately correct (PAC) learning}.
		
		\item A combinatorial tameness notion from Model Theory, which currently attracts considerable research activity, is known as \emph{`not the independence property' (NIP)}.\footnote{See Simon~\cite{Simon} for an introduction to NIP theories.} As o-minimality implies NIP, all results that hold for NIP structures are, in particular, applicable to any o-minimal structure. This observation is exploited in \cite[Proposition~4.5]{Krapp-Wirth} in order to ensure finite \emph{Vapnik--Chervonenkis (VC) dimensions} in the definable context. Finite VC dimensions, in turn, are the main requirement for the applicability of the Fundamental Theorem of Statistical Learning in order to derive PAC learnability.
	\end{enumerate}
	
	While all notions above that stem from Statistical Learning Theory are not subject of this research note, the two observations prompt us to ask whether the combinatorial tameness property NIP is already enough to guarantee measurability. We make this precise in our following main question. \pagebreak
	
	\begin{question}\label{question:main}
		Let $K$ be an NIP ordered field and let $D\subseteq K$ be definable in $K$. Is $D$ necessarily a Borel set with respect to the order topology on $K$?
	\end{question}
	
	In Question~\ref{question:main}, both the notion of NIP and of definability are regarded in the language of rings. A positive answer to this question is necessary for a strengthening of the Fundamental Theorem of Statistical Learning in the definable context (as presented in \cite[Theorem~4.7]{Krapp-Wirth}) from o-minimal to NIP structures to be feasible.
	We also point out that Question~\ref{question:main} relates to Shelah's Conjecture on the classification of NIP fields (cf.\ e.g.\ Dupont, Hasson and Kuhlmann~\cite[page~820]{Dupont-Hasson-Kuhlmann} and Johnson~\cite[Conjecture~1.9]{Johnson}). This connection will be elaborated on further in Section~\ref{sec::further-work}.
	
	As an initial step towards an investigation of Question~\ref{question:main}, we present in Theorem~\ref{theorem::main-result}, as the main result of our research note, a subfield $K$ of $\mathbb{R}$ that defines (in the language of rings) a non-Borel set.
	The field $K$ is a purely transcendental extension of $\mathbb{Q}$ by continuum many algebraically independent elements.
	While the definability of a non-Borel set immediately yields that $K$ cannot be o-minimal, we also point out further tameness properties that $K$ fails to exhibit (see Remark~\ref{remark::orderings}). Most crucially, the field $K$ has the independence property (i.e.\ it is not NIP). Moreover, $K$ is undecidable, not almost real closed, and it admits $2^\mathfrak{c}$ many pairwise non-isomorphic archimedean orderings and $2^\mathfrak{c}$ many pairwise non-isomorphic non-archimedean orderings.\footnote{We denote by $\mathfrak{c}=2^{\aleph_0}$ the cardinality of the continuum.}  The construction of the field $K$ heavily relies on Proposition~\ref{proposition::technical-heart}, which is the technical heart of our work. This proposition gives rise to a large family of pairwise disjoint algebraically independent subsets of~$\mathbb{R}$, all of which are pathological with respect to Borel-measurability. 
	
	\section{Preparatory Results}\label{sec::preliminaries}
	
	We assume ZFC in order to ensure the applicability of set-theoretic concepts in the context of ordinals. In particular, we will define a non-Borel set, employing the principle of transfinite induction. 
	
	We denote by $\mathbb{N}$ the set of positive natural numbers and set $\mathbb{N}_0=\mathbb{N}\cup\{0\}$. We use standard notations for intervals. For instance, given a set $X$ (linearly) ordered by $<$, and $a,b\in X$, we write $[a,b)$ for $\{x\in X\mid a\leq x< b\}$ and $X_{>a}$ for $\{x\in X\mid a<x\}$. We use the convention that $\emptyset$ and $X$ are also intervals.
	
	\subsection{Constructing a Suitable Family of Non-Borel Sets}
	
	In this section, we gather basic concepts and results from Measure Theory that underlie the proof of the main result of this note. The technical heart of this work is established in Proposition~\ref{proposition::technical-heart}. We assume some familiarity with $\sigma$--algebras and measures and refer the reader to Bogachev~\cite{Bogachev1} and Cohn~\cite{Cohn} for further details.
	
	Our measure-theoretic examination is set in the measure space $(\mathbb{R},\mathcal{B}(\mathbb{R}),\mu)$, where $\mathcal{B}(\mathbb{R})$ denotes the Borel $\sigma$--algebra on $\mathbb{R}$ and $\mu$ denotes the unique measure on $\mathcal{B}(\mathbb{R})$ that maps an interval $I=(a,b)$ with $a,b\in\mathbb{R}$ and $a<b$ to its length $b-a$ (cf.\ \cite[Corollary~1.5.9]{Bogachev1}). We denote by $(\mathbb{R}, \mathcal{L}(\mathbb{R}),\lambda)$ the Lebesgue completion of the measure space $(\mathbb{R}, \mathcal{B}(\mathbb{R}),\mu)$ as developed in \cite[\S\,1.5]{Bogachev1}. We recall that the Lebesgue measure $\lambda$ extends the Borel measure $\mu$. More precisely, we have $\mathcal{B}(\mathbb{R})\subseteq \mathcal{L}(\mathbb{R})$, and the restriction of $\lambda$ to $\mathcal{B}(\mathbb{R})$ coincides with $\mu$ (cf.\ [1, Theorem 1.5.6]). Furthermore, we recall $|\mathcal{B}(\mathbb{R})|=\mathfrak{c}<2^\mathfrak{c}=|\mathcal{L}(\mathbb{R})|$ (see Srivastava~\cite[page~104]{Srivastava}).
	
	Given a subset $M\subseteq\mathbb{R}$, we define the \textbf{inner measure} of $M$ by
	\begin{align*}
		\mu_*(M):=\max\{\mu(B)\mid B\in\mathcal{B}(\mathbb{R}), B\subseteq M\}\in[0,\infty],
	\end{align*}
	and the \textbf{outer measure} of $M$ by
	\begin{align*}
		\mu^*(M):=\min\{\mu(B)\mid B\in\mathcal{B}(\mathbb{R}),M\subseteq B\}\in[0,\infty],
	\end{align*}
	where $[0,\infty]=\mathbb{R}_{\geq 0}\cup\{\infty\}$. See \cite[page~38\,f.]{Cohn} for the definition in general measure spaces. Due to \cite[\S\,1.5, Exercise~5]{Cohn}, the inner and outer measure of a subset $M\subseteq \mathbb{R}$, as introduced above, are well-defined. Moreover, we note that $\mu_*(M)\leq\mu^*(M)$ for any $M\subseteq\mathbb{R}$, and that $\mu_*(M)=\mu^*(M)=\mu(M)$ for any $M\in\mathcal{B}(\mathbb{R})$. Thus, given a subset $M\subseteq\mathbb{R}$ with $\mu_*(M)<\mu^*(M)$, we immediately obtain $M\notin \mathcal{B}(\mathbb{R})$.  
	
	We now turn towards establishing Proposition~\ref{proposition::technical-heart}.
	
	
	\begin{lemma}\label{lemma::complement-family-cardinality}
		Let $A\subseteq\mathbb{R}$ be such that $\mu_*(\mathbb{R}\setminus A)>0$. Then the family $\{B\in\mathcal{B}(\mathbb{R})\mid A\subseteq B, \mu(\mathbb{R}\setminus B)>0\}$ has cardinality $\mathfrak{c}$.
	\end{lemma}
	\begin{proof}
		We choose $C\in\mathcal{B}(\mathbb{R})$ such that $C\subseteq (\mathbb{R}\setminus A)$ and $\mu(C)=\mu_*(\mathbb{R}\setminus A)$. By countable additivity of $\mu$, there exists $n\in\mathbb{Z}$ such that $\mu(C\cap[n,n+1])>0$. Consider the map
		\begin{align*}
			\phi\colon [n,n+1]&\to[0,\infty), \\
			x&\mapsto \mu(C\cap [n,x]).
		\end{align*}
		For $x,y\in[n,n+1]$ with $x\leq y$ we compute 
		$$0\leq \phi(y)-\phi(x)=\mu(C\cap [x,y])\leq\mu([x,y])=y-x.$$
		This shows that $\phi$ is $1$--Lipschitz and thus continuous. By the Intermediate Value Theorem, the image of $\phi$ is an interval of cardinality $\mathfrak{c}$, as $$\phi(n)=0<\mu(C\cap[n,n+1])=\phi(n+1).$$
		Thus, also the family $\mathcal{C}=\{C\cap[n,x]\mid x\in[n,n+1],\mu(C\cap[n,x])>0\}$ has cardinality $\mathfrak{c}$. Now, the fact that $\mathcal{C}\to \{B\in\mathcal{B}(\mathbb{R})\mid A\subseteq B, \mu(\mathbb{R}\setminus B)>0\},\linebreak
		M\mapsto \mathbb{R}\setminus M$
		is a well-defined injective map implies the claim.
	\end{proof}
	
	The following result due to Steinhaus~\cite[page~99\,f.]{Steinhaus} is commonly referred to as Steinhaus' Theorem (see Srivastava~\cite[Theorem~3.4.17]{Srivastava}).
	
	\begin{fact}[Steinhaus' Theorem]\label{fact::Steinhaus}
		Let $A\in\mathcal{L}(\mathbb{R})$ be such that $\lambda(A)>0$. Then the Minkowski difference $A-A=\{a-a'\mid a,a'\in A\}$ is a neighborhood of $0$.
	\end{fact}
	
	We now present the technical construction that the proof of our main result (Theorem~\ref{theorem::main-result}) relies on.
	
	\begin{proposition}\label{proposition::technical-heart}
		Let $K$ be a subfield of $\mathbb{R}$ such that the transcendence degree of $\mathbb{R}$ over $K$ is $\mathfrak{c}$. Then there exists a family $(A_\alpha)_{\alpha<\mathfrak{c}}$ of pairwise disjoint subsets of $\mathbb{R}$ satisfying the following conditions:
		\begin{enumerate}[(I)]
			\item\label{item::inner-zero} $\mu_*(A_\alpha)=0$ for any $\alpha<\mathfrak{c}$.
			\item\label{item::complement-inner-zero} $\mu_*(\mathbb{R} \setminus A_\alpha)=0$ for any $\alpha<\mathfrak{c}$.
			\item\label{item::algebraic-independence} The union $\dot{\bigcup\limits_{\alpha<\mathfrak{c}}}\,A_\alpha $ is algebraically independent over $K$.
		\end{enumerate}
	\end{proposition}
	\begin{proof}
		We consider the sets
		\begin{align*}
			\mathcal{B}_0&=\{B\in\mathcal{B}(\mathbb{R})\mid \mu(B)>0\}, \\
			\mathcal{B}_1&=\{C\in\mathcal{B}(\mathbb{R})\mid \mu(\mathbb{R}\setminus C)>0\}.
		\end{align*}
		As $|\mathcal{B}(\mathbb{R})|=\mathfrak{c}$ and $(0,x)\in \mathcal{B}_0\cap\mathcal{B}_1$ for any $x\in \mathbb{R}_{>0}$, we have $|\mathcal{B}_0|=|\mathcal{B}_1|=\mathfrak{c}$. Therefore, we can enumerate the sets and write $\mathcal{B}_0=\{B_\sigma\}_{\sigma<\mathfrak{c}}$ as well as $\mathcal{B}_1=\{C_\sigma\}_{\sigma<\mathfrak{c}}$.
		
		In the following, we proceed by transfinite induction to prove the existence of a single-indexed sequence $(w_\sigma)_{\sigma<\mathfrak{c}}$ and a double-indexed sequence $(a_{\sigma,\alpha})_{\alpha\leq\sigma<\mathfrak{c}}$ of real numbers. More precisely, there is an ``outer'' induction ranging over the index $\sigma$ and an ``inner'' induction ranging over the index~$\alpha$, which is bounded by $\sigma$. The index $\alpha$ indicates in which $A_\alpha$ the number $a_{\sigma,\alpha}$ will be included, i.e.\ we will later set $A_\alpha=\{a_{\sigma,\alpha}\mid \alpha\leq\sigma<\mathfrak{c}\}$. The auxiliary sequence $(w_\sigma)_{\sigma<\mathfrak{c}}$ will consist of numbers that are never included in the main sequence $(a_{\sigma,\alpha})_{\alpha\leq\sigma<\mathfrak{c}}$. The structure of this ``nested'' induction on $\sigma$ and $\alpha$ is visualized in Table~\ref{table::induction}.
		
		\begingroup
		\def\arraystretch{1.5}
		\centering
		\begin{table}[h!]
			\caption{Visualization of the Transfinite Induction.}\vspace{0.5em}
			\label{table::induction}
			\resizebox{\columnwidth}{!}{
				\begin{tabular}{|c||c||l|l|l|l|l|c|}
					\hline
					\diagbox{``outer'' induction $\sigma$}{``inner'' induction \\ \vspace{0.1em} $\alpha\leq \sigma$} & \cellcolor{black} &\multicolumn{1}{c|}{$0$} & \multicolumn{1}{c|}{$1$} & \multicolumn{1}{c|}{$2$} & \multicolumn{1}{c|}{$3$} & \multicolumn{1}{c|}{$4$} & \multicolumn{1}{c|}{$\cdots$} \\ \hline\hline
					
					\hspace{- 1.7em}$0$ & $w_0$ & $a_{0,0}$ & \cellcolor{black} & \cellcolor{black} & \cellcolor{black}  & \cellcolor{black}  & \cellcolor{black} \\ \hline
					
					\hspace{- 1.7em}$1$ & $w_1$ & $a_{1,0}$ & $a_{1,1}$ & \cellcolor{black} & \cellcolor{black}  & \cellcolor{black}  & \cellcolor{black} \\ \hline
					
					\hspace{- 1.7em}$2$ & $w_2$ & $a_{2,0}$ & $a_{2,1}$ & $a_{2,2}$ & \cellcolor{black}  & \cellcolor{black}  & \cellcolor{black} \\ \hline
					
					\hspace{- 1.7em}$3$ & $w_3$ & $a_{3,0}$ & $a_{3,1}$ & $a_{3,2}$ & $a_{3,3}$  & \cellcolor{black}  & \cellcolor{black} \\ \hline
					
					\hspace{- 1.7em}$4$ & $w_4$ & $a_{4,0}$ & $a_{4,1}$ & $a_{4,2}$ & $a_{4,3}$  & $a_{4,4}$  & \cellcolor{black} \\ \hline
					
					\hspace{- 1.7em}$\vdots$ & $\vdots$ & \multicolumn{1}{c|}{$\cdots$} & \multicolumn{1}{c|}{$\cdots$} &  \multicolumn{1}{c|}{$\cdots$} & \multicolumn{1}{c|}{$\cdots$}  & \multicolumn{1}{c|}{$\cdots$}  & \multicolumn{1}{c|}{$\cdots$} \\ \hline
					
					\multicolumn{1}{c}{} & \multicolumn{1}{c}{} & \multicolumn{1}{c}{$\uparrow$} & \multicolumn{1}{c}{$\uparrow$} & \multicolumn{1}{c}{$\uparrow$} & \multicolumn{1}{c}{$\uparrow$} & \multicolumn{1}{c}{$\uparrow$} &  \multicolumn{1}{c}{$\cdots$} \\
					
					\multicolumn{1}{c}{} & \multicolumn{1}{c}{} & \multicolumn{1}{c}{$A_0$} & \multicolumn{1}{c}{$A_1$} & \multicolumn{1}{c}{$A_2$} & \multicolumn{1}{c}{$A_3$} & \multicolumn{1}{c}{$A_4$} &  \multicolumn{1}{c}{$\cdots$} \\
			\end{tabular}}
			\vspace{0.7em}
			\centering
			\begin{minipage}[c]{0.8\textwidth}
				\small
				The table is filled row by row. For instance, $a_{3,1}$ is the $8^\mathrm{th}$ element of the main sequence  we fix. At the beginning of each step of the ``outer'' induction, we fix an element of the auxiliary sequence. At the end of each step of the ``outer'' induction, we begin to ``fill'' one new $A_\alpha$. More precisely, at the end of step $\sigma$ of the ``outer'' induction, we fix $a_{\sigma,\sigma}$, which will be a member of $A_\sigma$. 
			\end{minipage}
		\end{table}
		\endgroup
		
		The numbers in the sequence $(a_{\sigma,\alpha})_{\alpha\leq\sigma<\mathfrak{c}}$ are chosen carefully to ensure the following:
		\begin{itemize}
			\item $A_\alpha$ contains no element of $\mathcal{B}_0$ as a subset, which will imply \eqref{item::inner-zero}.
			
			\item The cardinality of the family $\{C\in\mathcal{B}_1\mid A_\alpha\subseteq C\}$ is strictly less than~$\mathfrak{c}$, which will imply \eqref{item::complement-inner-zero} (see Lemma~\ref{lemma::complement-family-cardinality}).
			
			\item Algebraic independence is maintained, which will imply \eqref{item::algebraic-independence}.
		\end{itemize}
		Given $\alpha\leq\sigma<\mathfrak{c}$, the part of the main sequence that we have already fixed when we are about to choose $a_{\sigma,\alpha}$ is given by the set
		$$D_{\sigma,\alpha}:=
		\{a_{\tau,\gamma}\mid \gamma\leq\tau<\sigma\}\cup\{a_{\sigma,\gamma}\mid \gamma<\alpha\}.$$
		Note that $D_{0,0}=\emptyset$. In the case $\sigma>0$, the first set $\{a_{\tau,\gamma}\mid \gamma\leq\tau<\sigma\}$ has been defined during earlier steps of the ``outer'' induction, while the second set $\{a_{\sigma,\gamma}\mid \gamma<\alpha\}$ has been defined at the current step of the ``outer induction'', but during earlier steps of the ``inner'' induction.
		
		Now, let $\sigma<\mathfrak{c}$, and assume that we have already fixed sequences $(w_\tau)_{\tau<\sigma}$ and $(a_{\tau,\gamma})_{\gamma\leq\tau<\sigma}$ satisfying the following conditions:
		\begin{enumerate}[(i)]
			\item\label{item::cond1} $w_\tau\in B_\tau$ for any $\tau<\sigma$,
			\item $a_{\tau,\gamma}\in (\mathbb{R}\setminus C_\tau)$ for any $\gamma\leq\tau<\sigma$,
			\item $a_{\tau,\gamma}\neq w_{\tau'}$ for any $\gamma\leq\tau<\sigma$ and $\tau'<\sigma$, 
			\item $a_{\tau,\gamma}\neq a_{\tau',\gamma'}$ for any $\gamma\leq\tau<\sigma$ and $\gamma'\leq \tau'<\sigma$ with $(\tau,\gamma)\neq(\tau',\gamma')$,
			\item\label{item::cond5} $\{a_{\tau,\gamma}\mid \gamma\leq\tau<\sigma\}$ is algebraically independent over $K$.
		\end{enumerate}
		We now choose suitable $w_\sigma$ and $(a_{\sigma,\alpha})_{\alpha\leq\sigma}$.
		\begin{itemize}
			\item Choice of $w_\sigma$: \\
			The set $B_\sigma\in \mathcal{B}_0$ has positive measure $\mu(B_\sigma)>0$, and thus we obtain $|B_\sigma|=\mathfrak{c}$ by applying Kechris~\cite[Theorem~13.6]{Kechris}. On the other hand, we compute 
			$$|D_{\sigma,0}|=|\{a_{\tau,\gamma}\mid \gamma\leq\tau<\sigma\}|\leq|\sigma|^2<\mathfrak{c}.$$
			Therefore, we can choose an element $w_\sigma$ from the set $B_\sigma\setminus D_{\sigma,0}$. 
			
			\item Choice of $(a_{\sigma,\alpha})_{\alpha\leq\sigma}$: \\
			Let $\alpha\leq\sigma$, and assume that we have already fixed a sequence $(a_{\sigma,\gamma})_{\gamma<\alpha}$ in an appropriate way (i.e.\ in accordance with suitable extensions of conditions \eqref{item::cond1}--\eqref{item::cond5}). We compute
			$$|D_{\sigma,\alpha}|=|\{a_{\tau,\gamma}\mid \gamma\leq\tau<\sigma\}|+|\{a_{\sigma,\gamma}\mid \gamma<\alpha\}|\leq |\sigma|^2+|\alpha|<\mathfrak{c}.$$
			Since the transcendence degree of $\mathbb{R}$ over $K$ is $\mathfrak{c}$, there exists a set $T_{\sigma,\alpha}\subseteq\mathbb{R}$ with $|T_{\sigma,\alpha}|=\mathfrak{c}$ such that $D_{\sigma,\alpha}\mathbin{\dot{\cup}}T_{\sigma,\alpha}$ is algebraically independent over $K$. The complement of the set $C_\sigma\in\mathcal{B}_1$ has positive measure $\mu(\mathbb{R}\setminus C_\sigma)>0$. Thus, by Fact~\ref{fact::Steinhaus} the Minkowski difference $(\mathbb{R}\setminus C_\sigma)-(\mathbb{R}\setminus C_\sigma)$ is a neighborhood of $0$. Since multiplication with non-zero elements of $K$ preserves algebraic independence over $K$, we can without loss of generality assume that $T_{\sigma,\alpha}$ is a subset of $(\mathbb{R}\setminus C_\sigma)-(\mathbb{R}\setminus C_\sigma)$. Hence, any $t\in T_{\sigma,\alpha}$ can be written as $t=u-v$ for some $u,v\in (\mathbb{R}\setminus C_\sigma)$. Set $U_t=\{x\in \{u,v\}\mid D_{\sigma,\alpha}\mathbin{\dot{\cup}} \{x\}\text{ is algebraically independent over } K\}$, and note that $U_t$ is non-empty, as $t=u-v$ and $D_{\sigma,\alpha}\mathbin{\dot{\cup}} \{t\}$ is algebraically independent over $K$. Now, consider the set $$U_{\sigma,\alpha}= \bigcup\limits_{t\in T_{\sigma,\alpha}}U_t.$$ Then the union $(D_{\sigma,\alpha}\mathbin{\dot{\cup}} U_{\sigma,\alpha})\cup\{t\}$ is algebraically dependent over $K$ for any $t\in T_{\sigma,\alpha}$, i.e.\
			$T_{\sigma,\alpha}$ is contained in the relative algebraic closure of the field $K(D_{\sigma,\alpha}\mathbin{\dot{\cup}} U_{\sigma,\alpha})\subseteq \mathbb{R}$. Since $|D_{\sigma,\alpha}|<|T_{\sigma,\alpha}|=\mathfrak{c}$, this yields
			$|U_{\sigma,\alpha}|=\mathfrak{c}$. Moreover, we have $|\{w_\tau\mid \tau\leq\sigma\}|<\mathfrak{c}$. Therefore, we can choose an element $a_{\sigma,\alpha}$ from the set $U_{\sigma,\alpha}\setminus\{w_\tau\mid \tau\leq\sigma\}$.
		\end{itemize}
		
		For $\alpha<\mathfrak{c}$ we now define 
		$$A_\alpha:=\{a_{\sigma,\alpha}\mid \alpha\leq\sigma<\mathfrak{c}\},$$
		and we verify the conditions \eqref{item::inner-zero}, \eqref{item::complement-inner-zero} and \eqref{item::algebraic-independence}:
		\begin{enumerate}[(I)]
			\item Let $\alpha<\mathfrak{c}$. To derive $\mu_*(A_\alpha)=0$, it suffices to verify $w_\sigma\in B_\sigma\setminus A_\alpha$ for any $\sigma<\mathfrak{c}$. Indeed, this implies $B\not\subseteq A_\alpha$ for any set $B\in\mathcal{B}_0$. Thus, let $\sigma<\mathfrak{c}$. We have $w_\sigma\in B_\sigma$ by our choice of $w_\sigma$. To verify $w_\sigma\notin A_\alpha$ we have to show that $w_\sigma\neq a_{\tau,\alpha}$ for any $\alpha\leq\tau<\mathfrak{c}$. Thus, let $\alpha\leq\tau<\mathfrak{c}$. If $\tau<\sigma$, then our choice of $w_\sigma$ ensures that $w_\sigma\neq a_{\tau,\alpha}$, as $w_\sigma\notin D_{\sigma,0}$ but $a_{\tau,\alpha}\in D_{\sigma,0}$. If $\sigma\leq\tau$, then our choice of $a_{\tau,\alpha}$ ensures that $w_\sigma\neq a_{\tau,\alpha}$, as 
			$a_{\tau,\alpha}\notin\{w_{\tau'}\mid \tau'\leq \tau\}$.
			
			\item Let $\alpha\leq\sigma<\mathfrak{c}$. Then $a_{\sigma,\alpha}\in A_\alpha\setminus C_\sigma$, and hence $A_\alpha\not\subseteq C_\sigma$. As a consequence, we obtain
			\begin{align*}
				&\,\{C\in\mathcal{B}_1\mid A_\alpha\subseteq C\}\\
				=&\,\{C_\sigma\mid \sigma<\mathfrak{c}, A_\alpha\subseteq C_\sigma\} \\
				\subseteq&\,\{C_\sigma\mid \sigma<\alpha\},
			\end{align*}
			which implies $|\{C\in\mathcal{B}_1\mid A_\alpha\subseteq C\}|\leq|\alpha|<\mathfrak{c}$. Applying Lemma~\ref{lemma::complement-family-cardinality} therefore yields $\mu_*(\mathbb{R}\setminus A_\alpha)=0$.
			
			\item Our choice of the sequence ensures that the set $D_{\sigma,\alpha}\cup\{a_{\sigma,\alpha}\}$ is algebraically independent over $K$ for any $\alpha\leq\sigma<\mathfrak{c}$, where $a_{\sigma,\alpha}$ is distinct from all elements of $D_{\sigma,\alpha}$. Thus, the members of the sequence $(A_\alpha)_{\alpha<{\mathfrak{c}}}$ are pairwise disjoint, and we obtain the algebraic independence of the union $\dot{\bigcup\limits_{\alpha<\mathfrak{c}}}\,A_\alpha$
			over $K$. \qedhere
		\end{enumerate}
	\end{proof}
	
	\pagebreak
	\begin{remark}\label{remark::heart}\
		\begin{enumerate}[(a)]
			\item We note that Proposition~\ref{proposition::technical-heart} applies to any subfield $K$ of $\mathbb{R}$ with $|K|<\mathfrak{c}$, such as $\mathbb{Q}$ or the relative algebraic closure of $\mathbb{Q}$ in $\mathbb{R}$.
			
			\item For each $\alpha<\mathfrak{c}$, the set $A_\alpha$ is non-Borel, i.e.\ $A_\alpha\notin\mathcal{B}(\mathbb{R})$. Indeed, if $A_\alpha$ were Borel, then $A_\alpha, (\mathbb{R}\setminus A_\alpha)\in\mathcal{B}(\mathbb{R})$ yields the contradiction $$\mu(\mathbb{R})=\mu(A_\alpha)+\mu(\mathbb{R}\setminus A_\alpha)=\mu_*(A_\alpha)+\mu_*(\mathbb{R}\setminus A_\alpha)=0.$$
		
		\item\label{item::cardinal} By way of construction, for each $\alpha<\mathfrak{c}$, the cardinality of $A_\alpha$ is given by $|\mathfrak{c}\setminus \alpha|=\mathfrak{c}$.
	\end{enumerate}
\end{remark}

We complete this section with a lemma that will be applied in the proof of Theorem~\ref{theorem::main-result}.

\begin{lemma}\label{lemma::inner-measure-zero}
	Let $A\subseteq\mathbb{R}$. Then the following conditions are equivalent:
	\begin{enumerate}[(i)]
		\item\label{item::positive-inner} $\mu_*(A)=0$.
		\item\label{item::outer-intervals} $\mu^*(I\setminus A)=\mu(I)$ for any interval $I$ of $\mathbb{R}$.
	\end{enumerate}
\end{lemma}
\begin{proof}
	We prove both directions by contraposition. If $\mu_*(A)>0$, then there exists $B\in\mathcal{B}(\mathbb{R})$ such that $B\subseteq A$ and $\mu(B)>0$. Recall that $\mathbb{R}$ can be partitioned into countably many intervals of finite length. Thus, by countable additivity of $\mu$, there must exist an interval $I\subseteq\mathbb{R}$ of finite length such that $0<\mu(I\cap B)$. Since $(I\setminus A)\subseteq (I\setminus B)$ and $(I\setminus B)\in\mathcal{B}(\mathbb{R})$, we obtain 
	$$\mu^*(I\setminus A)\leq \mu(I\setminus B)=\mu(I)-\mu(I\cap B)<\mu(I)$$
	(cf.\ Cohn~\cite[Proposition~1.2.1]{Cohn}). For the converse implication, let $I\subseteq \mathbb{R}$ be an interval with $\mu^*(I\setminus A)<\mu(I)$. Further, we choose $B\in\mathcal{B}(\mathbb{R})$ such that $(I\setminus A)\subseteq B$ and $\mu(B)=\mu^*(I\setminus A)$. Since we can write $I=(I\cap B)\mathbin{\dot{\cup}}(I\setminus B)$ and we have 
	$$\mu(I\cap B)\leq \mu(B)=\mu^*(I\setminus A)<\mu(I),$$
	the additivity of $\mu$ implies $\mu(I\setminus B)>0$. This yields $\mu_*(A)>0$, as $(I\setminus B)\subseteq A$ and $(I\setminus B)\in \mathcal{B}(\mathbb{R})$.
\end{proof}

\subsection{First-Order Defining the Integers}

In this section, we briefly introduce the general model-theoretic setup and terminology. We assume some familiarity with first-order logic and refer the reader to Marker~\cite{Marker} and Poizat~\cite{Poizat} for further details.

Our model-theoretic notions are all set in the \textbf{language of rings}, which is given by $\mathcal{L}_\mathrm{r}=\{+,-,\cdot,0,1\}$. Whenever the interpretation of the symbols is clear from the context, we simply write the domain instead of the structure, i.e.\ the $\mathcal{L}_\mathrm{r}$--structure of a field $(K,+,-,\cdot,0,1)$ is simply written as $K$. 
Given a field $K$ and $n\in \mathbb{N}$, 
a set $A\subseteq K^n$ is called \textbf{definable} if there is an $\mathcal{L}_\mathrm{r}$--formula (potentially with parameters from $K$) that defines $A$ over $K$. We point out in particular when definability is obtained without parameters and then use the term \textbf{$\emptyset$--definable}.

In the following, we only introduce the model-theoretic concept \emph{IP}, which is originally due to Shelah~\cite{Shelah1971}, for fields. See Poizat~\cite[\S\,12.4]{Poizat} for a general definition.
Let $K$ be a field and let $\varphi(x_1,\ldots,x_n;y_1,\ldots,y_\ell)$ be an $\mathcal{L}_\mathrm{r}$--formula. Then $\varphi(\underline{x};\underline{y})$ has the \textbf{independence property (IP)} over $K$ if for any $m\in\mathbb{N}$ there is a set $\{\underline{a}_1,\ldots,\underline{a}_m\}\subseteq K^n$ and a set $\{\underline{b}_I\mid I\subseteq \{1,\dots,m\}\}\subseteq K^\ell$ such that for any $J\subseteq \{1,\dots,m\}$ we have
$$K\models \bigwedge_{i\in J} \varphi(\underline{a}_i;\underline{b}_{J}) \wedge \bigwedge_{i\in \{1,\dots,m\}\setminus J} \neg \varphi(\underline{a}_i;\underline{b}_{J}).$$
The field $K$ has the \textbf{independence property} if there is an $\mathcal{L}_\mathrm{r}$--formula that has the independence property over $K$. Moreover, $K$ is said to be \textbf{NIP} if it does not have the independence property. 

We provide in Corollary~\ref{corollary:trans-IP} the sufficient condition we use in Theorem~\ref{theorem::main-result} to verify the independence property.

\begin{observation}\label{oberservation::IP}
	Let $K$ be a field of characteristic $0$ in which $\mathbb{Z}$ is $\emptyset$--definable.\footnote{This observation also holds when $\mathbb{Z}$ is definable with parameters. However, we point out $\emptyset$--definability here, as our main definability results do not need any parameters.} Then $K$ has the independence property.
\end{observation}
\begin{proof}
	Let $\psi(x)$ be an $\mathcal{L}_\mathrm{r}$--formula defining $\mathbb{Z}$ in $K$. The $\mathcal{L}_\mathrm{r}$--structure $\mathbb{Z}$ has the independence property, as witnessed by the $\mathcal{L}_\mathrm{r}$--formula $\exists z\, x=y\cdot z$ (see e.g.\ Poizat~\cite[\S\,12.4]{Poizat}). Hence, the $\mathcal{L}_\mathrm{r}$--formula $\varphi(x;y)$ given by $\exists z\, (\psi(z)\land x=y\cdot z)$ has the independence property over $K$.
\end{proof}

Recall that a field is \textbf{(formally) real} if it admits an ordering making it an ordered field. 
In order to obtain the independence property for the field in Theorem~\ref{theorem::main-result}, we combine Observation~\ref{oberservation::IP} with the following.

\begin{fact}[{{R.~Robinson~\cite[\S\,5]{R-Robinson}}}]\label{fact::robinson}
	Let $F$ be a real field that admits an archimedean ordering and let $t$ be transcendental over $F$. Then $\mathbb{Z}$ is $\emptyset$--definable in $F(t)$.
\end{fact}

\begin{corollary}\label{corollary:trans-IP}
	Let $F$ be a real field that admits an archimedean ordering and let $t$ be transcendental over $F$. Then $F(t)$ has the independence property.
\end{corollary}

Note that due to Hölder's Theorem\footnote{Applied to ordered fields, this theorem states that for any archimedean ordered field $(F,<_F)$, there is a \emph{unique} order-preserving $\mathcal{L}_{\mathrm{r}}$--embedding from $(F,<_F)$ into $(\mathbb{R},<)$. See Hölder~\cite[Erster Theil]{Hoelder}, Engler and Prestel~\cite[Proposition~2.1.1]{Engler-Prestel}.}, real fields that admit at least one archimedean ordering are precisely the subfields of $\mathbb{R}$ (up to $\mathcal{L}_{\mathrm{r}}$--i\-so\-mor\-phism).

\section{Building a Wild Ordered Subfield of the Reals}\label{sec::main}

In this section, we assume some familiarity with Borel $\sigma$--algebras in general topological spaces and refer the reader to Bogachev~\cite[\S\,6.2]{Bogachev2} for further details.

We now turn to the construction of a subfield $K$ of $\mathbb{R}$ that has the independence property and defines a non-Borel set $D\subseteq K$. In this context, given a subset $X\subseteq \mathbb{R}$, there are two natural candidates for Borel $\sigma$--algebras on $X$: 
\begin{itemize}
	\item We can endow $X$ with the \textbf{trace $\sigma$--algebra} 
	$$\mathcal{B}(X):=\{B\cap X\mid B\in\mathcal{B}(\mathbb{R})\}$$
	induced by $\mathcal{B}(\mathbb{R})$. Denoting by $\tau_\mathbb{R}$ the order topology\footnote{Recall that the order topology on $\mathbb{R}$ coincides with the Euclidean topology, and thus it generates $\mathcal{B}(\mathbb{R})$.} on $\mathbb{R}$, it follows from \cite[Lemma~6.2.4]{Bogachev2} that $\mathcal{B}(X)$ is precisely the Borel $\sigma$--algebra generated by the subspace topology $\{U\cap X\mid U\in\tau_\mathbb{R}\}$ on $X$ induced by~$\tau_\mathbb{R}$.
	
	\item The subset $X$ of $\mathbb{R}$ inherits the linear ordering of $\mathbb{R}$. Thus, $X$ is naturally endowed with the order topology $\tau_X$, and we can consider the Borel $\sigma$--algebra $\mathcal{B}(\tau_X)$ generated by $\tau_X$.
\end{itemize}

The following result shows that the two Borel $\sigma$--algebras considered above coincide.

\begin{lemma}\label{lemma::Borel-identities}
	For any subset $X\subseteq\mathbb{R}$, we have $\mathcal{B}(X)=\mathcal{B}(\tau_X)$.
\end{lemma}
\begin{proof}
	The inclusion $\mathcal{B}(\tau_X)\subseteq\mathcal{B}(X)$ follows from the fact that the generators of the $\sigma$--algebras satisfy $\tau_X\subseteq \{U\cap X\mid U\in\tau_\mathbb{R}\}$. For the proof of the other inclusion, we consider the set $$\Sigma:=\{B\in\mathcal{B}(\mathbb{R})\mid (B\cap X)\in\mathcal{B}(\tau_X)\}.$$
	Note that $\Sigma$ is a $\sigma$--algebra on $\mathbb{R}$. If we prove the inclusion $\mathcal{B}(\mathbb{R})\subseteq\Sigma$, then the claim of the lemma follows. Hence, it suffices to verify $\mathbb{R}_{\leq a}\in \Sigma$ for any $a\in\mathbb{R}$ (cf.\ Bogachev~\cite[Lemma~1.2.11]{Bogachev1}). Given $a\in\mathbb{R}$, set 
	$b=\sup (\mathbb{R}_{\leq a}\cap X),$ 
	where the supremum is taken in $\mathbb{R}$. If $b\in X$, then $(\mathbb{R}_{\leq a}\cap X)=X_{\leq b}\in\mathcal{B}(\tau_K)$. Otherwise, for any $n\in\mathbb{N}$ there exists $b_n\in(\mathbb{R}_{\leq a}\cap X)$ such that $b-\frac{1}{n}<b_n< b\leq a$. Hence, 
	$$\mathbb{R}_{\leq a}\cap X=\bigcup\limits_{n\in\mathbb{N}} X_{\leq b_n}\in\mathcal{B}(\tau_{X}).$$
	This yields $\mathbb{R}_{\leq a}\in\Sigma$, as required.
\end{proof}

We now turn towards the proof of our main result. First, we establish a lemma, whose proof is straightforward basic algebra and therefore left to the reader.

\begin{lemma}\label{lemma::squares}
	Let $C,C'\subseteq\mathbb{R}$ be disjoint subsets such that $C\subseteq \mathbb{R}_{\geq 0}$ and $C\mathbin{\dot{\cup}} C'$ is algebraically independent over $\mathbb{Q}$, and let $K=\mathbb{Q}(\sqrt{C}\cup C')$, where $\sqrt{C}=\{\sqrt{c}\mid c\in C\}$. Then the following hold: 
	\begin{enumerate}[(i)]
		\item\label{item::square} Any $c'\in C'$ is not a square in $K$, i.e.\ $b^2\neq c'$ for any $b\in K$.
		\item\label{item::transcendental} The set $\sqrt{C}\mathbin{\dot{\cup}}C'$ is algebraically independent over $\mathbb{Q}$.
	\end{enumerate}
\end{lemma}

\begin{theorem}\label{theorem::main-result}
	There exists a subfield $K\subseteq \mathbb{R}$ of cardinality $\mathfrak{c}$ that has the independence property in the language of rings and $\emptyset$--defines a set $D\subseteq K$ with $D\notin\mathcal{B}(\tau_K)$.
\end{theorem}
\begin{proof}
	By Proposition~\ref{proposition::technical-heart} there exist two disjoint sets $A,A'\subseteq\mathbb{R}$ of cardinality $\mathfrak{c}$ (see Remark~\ref{remark::heart}\eqref{item::cardinal}) such that $\mu_*(A)=\mu_*(A')=\mu_*(\mathbb{R}\setminus A)=\mu_*(\mathbb{R}\setminus A')=0$ and $A\mathbin{\dot{\cup}} A'$ is algebraically independent over $\mathbb{Q}$. Set $\sqrt{A_{\geq 0}}=\{\sqrt{a}\mid a\in A_{\geq 0}\}$, and let
	$$K=\mathbb{Q}\!\left(\sqrt{A_{\geq 0}}\cup A'\right).$$ 
	For any $a'\in A'$, we can write $K=F(a')$, where $F=\mathbb{Q}(\sqrt{A_{\geq 0}}\cup (A'\setminus\{a'\}))$. By Lemma~\ref{lemma::squares}\eqref{item::transcendental}, $a'$ is transcendental over $F$. Corollary~\ref{corollary:trans-IP} yields that $K$ has the independence property. 
	
	Now, consider the set 
	$$D=\{y^2\mid y\in K\}\subseteq K$$
	of squares in $K$. Then the $\mathcal{L}_\mathrm{r}$--formula $\exists y\; x=y^2$ defines $D$ in $K$. By the definition of $K$, we have $\sqrt{A_{\geq 0}}\subseteq K$, which implies $A_{\geq 0}\subseteq D$. On the other hand, Lemma~\ref{lemma::squares}\eqref{item::square} implies $A'\cap D=\emptyset$, and thus $A'\subseteq (K\setminus D)$. In the following, we assume for a contradiction that $D\in \mathcal{B}(\tau_K)$. Then by Lemma~\ref{lemma::Borel-identities} there exists a Borel set $B\in\mathcal{B}(\mathbb{R})$ such that $D=B\cap K$. We compute $\mu(B\cap[0,1])$ as well as $\mu([0,1]\setminus B)$:
	\begin{itemize}
		\item 
		Since $\mu_*(\mathbb{R}\setminus A)=0$, applying Lemma~\ref{lemma::inner-measure-zero} yields $\mu^*(A\cap[0,1])=\mu([0,1])=1$. Moreover, $A_{\geq 0}\subseteq D\subseteq B$ implies $(A\cap[0,1])\subseteq (B\cap[0,1])$. We compute
		\begin{align*}
			1&=\mu^*(A\cap [0,1]) \\
			&=\min\{\mu(C)\mid C\in\mathcal{B}(\mathbb{R}), (A\cap[0,1])\subseteq C\} \\
			&\leq\mu(B\cap[0,1])\\
			&\leq\mu([0,1])=1,
		\end{align*}
		and hence obtain $\mu(B\cap[0,1])=1$.
		
		\item 
		Since $\mu_*(\mathbb{R}\setminus A')=0$, again $\mu^*(A'\cap[0,1])=\mu([0,1])=1$ by Lemma~\ref{lemma::inner-measure-zero}. Moreover, $A'\subseteq (K\setminus D)= (K\setminus B)\subseteq (\mathbb{R}\setminus B)$ implies $(A'\cap[0,1])\subseteq ([0,1]\setminus B)$. We compute
		\begin{align*}
			1&=\mu^*(A'\cap [0,1]) \\
			&=\min\{\mu(C)\mid C\in\mathcal{B}(\mathbb{R}), (A'\cap[0,1])\subseteq C\} \\
			&\leq\mu([0,1]\setminus B) \\
			&\leq\mu([0,1])=1,
		\end{align*}
		and hence obtain $\mu([0,1]\setminus B)=1$.
	\end{itemize}
	We obtain
	$$1=\mu([0,1])=\mu([0,1]\cap B)+\mu([0,1]\setminus B)=1+1=2,$$
	the required contradiction.
\end{proof}

We conclude this section with a final remark, in which we also collect further observations regarding tameness and related properties that the field $K$ from Theorem~\ref{theorem::main-result} fails to exhibit.

\begin{remark}\label{remark::orderings}\ 
	\begin{enumerate}[(a)]
		
		\item\label{item::cardinality-argument} Due to Lemma~\ref{lemma::squares}\eqref{item::transcendental}, as an $\mathcal{L}_{\mathrm{r}}$--structure the field $K$ is simply an extension of $\mathbb{Q}$ by continuum many algebraically independent elements $(\alpha_i)_{i<\mathfrak{c}}$ of $\mathbb{R}$. Thus, $K$ is $\mathcal{L}_{\mathrm{r}}$--isomorphic to the rational function field $\mathbb{Q}(t_i\mid i<\mathfrak{c})$ over $\mathbb{Q}$ in continuum many variables. We prove in Theorem~\ref{theorem::main-result} that the $\emptyset$--definable set $D$ of squares in $K$ is non-Borel with respect to the order topology, i.e.\ $D\notin\mathcal{B}(\tau_K)$. This immediately yields that $D$ is non-Borel in $\mathbb{R}$, i.e.\ $D\notin\mathcal{B}(\mathbb{R})$, since $\{B\in\mathcal{B}(\mathbb{R})\mid B\subseteq K\}\subseteq\mathcal{B}(\tau_K)$ (see Lemma~\ref{lemma::Borel-identities}). The existence of such a subfield $K\subseteq\mathbb{R}$, for which the set $D$ of squares in $K$ is not a member of $\mathcal{B}(\mathbb{R})$, can alternatively be established via the following cardinality argument: We obtain a set $\Phi$ of $2^\mathfrak{c}$ many pairwise distinct embeddings from $\mathbb{Q}(t_i\mid i<\mathfrak{c})$ into $\mathbb{R}$ by mapping each $t_i$ either to $\alpha_i$ or to~$\alpha_i^2$. For any $\varphi\in \Phi$, denote by $K_\varphi=\mathrm{Im}(\varphi)$ the obtained subfield of $\mathbb{R}$ and by $D_\varphi=\{y^2\mid y\in K_\varphi\}$ its set of squares. Depending on the choice of $\varphi$, for each $i<\mathfrak{c}$, we either have $\alpha_i^2\in D_\varphi$ or $\alpha_i^2\in K_\varphi\setminus D_\varphi$. Therefore, the $2^\mathfrak{c}$ many embeddings give rise to $2^\mathfrak{c}$ many subfields of~$\mathbb{R}$ with pairwise distinct sets of squares. On the other hand, we recall $|\mathcal{B}(\mathbb{R})|=\mathfrak{c}<2^{\mathfrak{c}}$. Hence, for $2^\mathfrak{c}$ many embeddings $\varphi$, we obtain a subfield $K_\varphi$ for which $D_\varphi\notin\mathcal{B}(\mathbb{R})$. However, this cardinality argument is not sufficient to establish the existence of a subfield $K\subseteq \mathbb{R}$ defining a set that is not Lebesgue measurable, since the Lebesgue $\sigma$--algebra $\mathcal{L}(\mathbb{R})$ has cardinality $2^\mathfrak{c}$. In contrast, our argumentation in the proof of Theorem~\ref{theorem::main-result} also applies to the Lebesgue measure $\lambda$ on $\mathbb{R}$. Since $(\mathbb{R}, \mathcal{L}(\mathbb{R}),\lambda)$ is the completion of the measure space $(\mathbb{R}, \mathcal{B}(\mathbb{R}),\mu)$ (cf.\ Bogachev~\cite[\S\,1.5]{Bogachev1}), it can be shown that the induced inner and outer measures coincide, i.e.\ $\mu_*=\lambda_*$ and $\mu^*=\lambda^*$. Therefore, the proof of Theorem~\ref{theorem::main-result} can easily be enhanced to yield $D\notin\{M\cap K\mid M\in\mathcal{L}(\mathbb{R})\}$, and thus $D\notin\mathcal{L}(\mathbb{R})$. Since we were mainly interested in Borel measurability, we did not treat the Lebesgue measure in this note.
		
		\item 
		In the following we outline how $K$ can be endowed with $2^{\mathfrak{c}}$ pairwise non-isomorphic\footnote{Two orderings $<$ and $<'$ on a real field $F$ are \textbf{isomorphic} if there exists an order-preserving $\mathcal{L}_{\mathrm{r}}$--isomorphism from $(F,<)$ to $(F,<')$.} archimedean and  $2^{\mathfrak{c}}$ pairwise non-isomorphic non-ar\-chi\-me\-dean orderings. Note that these cardinalities are maximal, as $K$ admits exactly $2^{\mathfrak{c}}$ binary relations. Consider the $2^{\mathfrak{c}}$ many distinct subfields $K_{\varphi}$ of $\mathbb{R}$ from \eqref{item::cardinality-argument} above. For each $\varphi\in\Phi$, as $K$ is $\mathcal{L}_{\mathrm{r}}$--isomorphic to  $K_{\varphi}$, we can endow $K$ with an ordering $<_\varphi$ via this isomorphism. This yields an order-preserving $\mathcal{L}_\mathrm{r}$--isomorphism from $(K_{\varphi},<)$ to $(K,<_\varphi)$, where $<$ denotes the ordering induced by the one on $\mathbb{R}$. Thus, for two distinct embeddings $\varphi,\varphi'\in\Phi$, the orderings $<_\varphi$ and $<_{\varphi'}$ are non-isomorphic due to Hölder's Theorem. Regarding the non-archimedean orderings, we note that the field $K$ is $\mathcal{L}_\mathrm{r}$--isomorphic to the rational function field $K(t)$ in one variable $t$, and each ordering on $K$ gives rise to a non-archimedean ordering on $K(t)$ by setting $t>K$. Thus, $K$ admits $2^{\mathfrak{c}}$ many pairwise non-isomorphic non-archimedean orderings.
		
		\item The field $K$ fails to exhibit certain further tameness properties. It cannot be o-minimal, since definable sets in o-minimal ordered fields are Borel with respect to the order topology (cf.\ Kaiser~\cite[Proposition~1.1]{Kaiser}, Karpinski and Macintyre~\cite[Lemma~6]{Karp-Mac}). In particular, $K$ is not real closed (see Pillay and Steinhorn~\cite[Proposition~1.4]{Pillay-Steinhorn}). This further implies that $K$ is not almost real closed. Indeed, any almost real closed field that admits at least one archimedean ordering is real closed (see Section~\ref{sec::further-work}).
		Moreover, the $\emptyset$--definability of $\mathbb{Z}$ implies that $K$ is undecidable (cf.\ J.~Robinson~\cite{J-Robinson}). 
	\end{enumerate}
\end{remark}

\section{Further Work}\label{sec::further-work}

From the perspective of Tame Geometry, the ordered field $K$ from Theorem~\ref{theorem::main-result} (with the ordering inherited from $\mathbb{R}$) is a pathological example of a ``wild'' field. Generally in Tame Geometry, dividing lines are of interest: Are there structures that exhibit one certain tameness property but fail to exhibit another certain tameness property? In light of this, we re-evaluate our main question from the introduction (Question~\ref{question:main}) for subfields of the reals. 

\begin{question}\label{question::NIP-nonBorel}
	Is there an NIP subfield $K$ of $\mathbb{R}$ that defines a set $D\subseteq K$ with $D\notin \mathcal{B}(K)$?
\end{question} 

Recall that if ``NIP'' is replaced by ``o-minimal'' in  Question~\ref{question::NIP-nonBorel}, then the answer to the resulting question is negative.

Next, we elaborate on how Question~\ref{question::NIP-nonBorel} relates to Shelah's Conjecture on the classification of NIP fields, as stated (up to minor variations) in Dupont, Hasson and Kuhlmann~\cite[page~820]{Dupont-Hasson-Kuhlmann} and Johnson~\cite[Conjecture~1.9]{Johnson}).
As we are interested in real fields, we record in Conjecture~\ref{conjecture::Shelah-real} below a suitable specialization of this conjecture (cf.\ Krapp, Kuhlmann and Lehéricy~\cite[Conjecture~6.2]{Krapp-Kuhlmann-Lehericy2} and \cite[Conjecture~1.2]{Krapp-Kuhlmann-Lehericy}) in terms of almost real closed fields.\footnote{In \cite{Krapp-Kuhlmann-Lehericy2} and \cite{Krapp-Kuhlmann-Lehericy} the conjectures are stated for \emph{strongly} NIP \emph{ordered} fields (in the language of ordered rings) rather than NIP real fields (in the language of rings). However, several accounts of Shelah's Conjecture are stated for NIP fields, such as \cite[Conjecture~1.9]{Johnson}. This conjecture can be specialized to NIP real fields in a similar fashion as in \cite{Krapp-Kuhlmann-Lehericy2} and \cite{Krapp-Kuhlmann-Lehericy}.} A real field $K$ is called \textbf{almost real closed} if it admits a henselian valuation with real closed residue field (see Delon and Farré~\cite{Delon-Farre}).

\begin{conjecture}\label{conjecture::Shelah-real}
	Any NIP real field is almost real closed.
\end{conjecture}

In light of Conjecture~\ref{conjecture::Shelah-real} and Question~\ref{question::NIP-nonBorel}, we ask the following question, to which a negative answer would be expected.

\begin{question}\label{question::ARC-nonBorel}
	Is there an almost real closed field $K$ that defines a set $D\subseteq K$ with $D\notin\mathcal{B}(\tau_K)$?
\end{question}

Here $\tau_K$ denotes the order topology with respect to any ordering on $K$. By Knebusch and Wright~\cite[Lemma~2.1]{Knebusch-Wright}, any henselian valuation on a real field~$K$ is convex with respect to any ordering on $K$. Therefore, if an almost real closed field $K$ admits at least one archimedean ordering, it follows that any henselian valuation on $K$ must be trivial, yielding that $K$ is real closed. Consequently, Conjecture~\ref{conjecture::Shelah-real} can be further specialized as follows:

\begin{conjecture}\label{conjecture::Shelah-arch}
	Let $K$ be an NIP real field that admits at least one ar\-chi\-me\-de\-an ordering. Then $K$ is real closed.
\end{conjecture}


As an alternative to Question~\ref{question::NIP-nonBorel}, one may consider the significance of the orderings the field admits:

\begin{question}\label{question::orderings-nonBorel}
	Is there a subfield $K$ of $\mathbb{R}$ that only admits archimedean orderings and that defines a set $D\subseteq K$ with $D\notin \mathcal{B}(K)$?
\end{question}

The respective first conditions on $K$ in Question~\ref{question::NIP-nonBorel}, Question~\ref{question::ARC-nonBorel} and Question~\ref{question::orderings-nonBorel} encode tameness requirements, which the field from Theorem~\ref{theorem::main-result} fails to meet. 


\begingroup
\renewcommand{\bibname}{References}

\endgroup

\pagebreak

\textbf{CRediT Authorship Contribution Statement:} Lothar Sebastian\linebreak Krapp: conceptualization (equal); funding acquisition (lead); investigation (supporting); project administration (supporting); supervision (lead); writing – original draft (supporting); writing – review \& editing (equal). \linebreak Matthieu Vermeil: conceptualization (equal); investigation (lead); writing – original draft (equal); writing – review \& editing (equal). Laura Wirth: conceptualization (equal); investigation (supporting); project administration (lead); writing – original draft (equal); writing – review \& editing (equal). 

\textbf{Funding:} The first and third author received partial project funding from Vector Stiftung as well as from the Network Platform \textit{Connecting Statistical Logic, Dynamical Systems and Optimization} of Universität Konstanz.

\textbf{Acknowledgments:} We started this research project during the \emph{Tame Geometry} workshop, as part of the thematic month titled \emph{Singularities, Differential Equations, Transcendence}, at CIRM in February 2025. We wish to thank the organizers of the workshop as well as CIRM for its hospitality. Moreover, we thank Philip Dittmann for providing an argument that led to Remark~\ref{remark::orderings}\eqref{item::cardinality-argument}. The third author is grateful to Carl Eggen for helpful comments on the presentation of the proof of Proposition~\ref{proposition::technical-heart}, to Lasse Vogel for valuable discussions on henselian valuations, and she would like to extend special thanks to Salma Kuhlmann for the supervision of her doctoral research project, which this work is part of.

\textbf{Conflict of Interest:} The authors declare no conflict of interest. The funders had no role in the design and conduct of the study; preparation, review, or approval of the manuscript; and decision to submit the manuscript for publication.

\textbf{Data and Materials Availability:} Not applicable.

\textbf{Code Availability:} Not applicable.

\textbf{Ethical Approval:} Not applicable.

\textbf{Consent to Participate:} Not applicable.

\textbf{Consent for Publication:} Not applicable. The research note does not include data or images that require permission to be published.

\end{document}